\documentclass[reqno,11pt]{amsart}
\usepackage{amscd,amssymb}
\usepackage{color}

\theoremstyle{plain}
\newtheorem{Thm}[subsection]{Theorem}
\newtheorem{Cor}[subsection]{Corollary}
\newtheorem{Lem}[subsection]{Lemma}
\newtheorem{Prop}[subsection]{Proposition}
\newtheorem{Conj}[subsection]{Conjecture}

\theoremstyle{definition}
\newtheorem{Def}[subsection]{Definition}

\theoremstyle{remark}

\newtheorem{Rem}[subsection]{Remark}

\newtheorem{conj}{Conjecture}

\numberwithin{equation}{section}
\renewcommand{\rm}{\normalshape}

\newif\ifShowLabels
\ShowLabelstrue
\newdimen\theight
\def\TeXref#1{%
    \leavevmode\vadjust{\setbox0=\hbox{{\tt
        \quad\quad  {\small \rm #1}}}%
    \theight=\ht0
    \advance\theight by \lineskip
    \kern -\theight \vbox to
    \theight{\rightline{\rlap{\box0}}%
    \vss}%
    }}%

\ShowLabelsfalse

\renewcommand{\sec}[2]{\section{#2}\label{S:#1}%
    \ifShowLabels \TeXref{{S:#1}} \fi}
\newcommand{\ssec}[2]{\subsection{#2}\label{SS:#1}%
    \ifShowLabels \TeXref{{SS:#1}} \fi}

\newcommand{\refs}[1]{Section ~\ref{S:#1}}
\newcommand{\refss}[1]{Section ~\ref{SS:#1}}

\newcommand{\reft}[1]{Theorem ~\ref{T:#1}}

\newcommand{\refp}[1]{Proposition ~\ref{P:#1}}
\newcommand{\refc}[1]{Corollary ~\ref{C:#1}}

\newcommand{\refe}[1]{\eqref{E:#1}}
\newcommand{\refco}[1]{Conjecture ~\ref{Co:#1}}

\newenvironment{thm}[1]%
    { \begin{Thm} \label{T:#1}  \ifShowLabels \TeXref{T:#1} \fi }%
    { \end{Thm} }

\renewcommand{\th}[1]{\begin{thm}{#1} \sl }
\renewcommand{\eth}{\end{thm} }

\newenvironment{lemma}[1]%
    { \begin{Lem} \label{L:#1}  \ifShowLabels \TeXref{L:#1} \fi }%
    { \end{Lem} }
\newcommand{\lem}[1]{\begin{lemma}{#1} \sl}
\newcommand{\elem}{\end{lemma}}

\newenvironment{propos}[1]%
    { \begin{Prop} \label{P:#1}  \ifShowLabels \TeXref{P:#1} \fi }%
    { \end{Prop} }
\newcommand{\prop}[1]{\begin{propos}{#1}\sl }
\newcommand{\eprop}{\end{propos}}

\newenvironment{corol}[1]%
    { \begin{Cor} \label{C:#1}  \ifShowLabels \TeXref{C:#1} \fi }%
    { \end{Cor} }
\newcommand{\cor}[1]{\begin{corol}{#1} \sl }
\newcommand{\ecor}{\end{corol}}

\newenvironment{defeni}[1]%
    { \begin{Def} \label{D:#1}  \ifShowLabels \TeXref{D:#1} \fi }%
    { \end{Def} }
\newcommand{\defe}[1]{\begin{defeni}{#1} \sl }
\newcommand{\edefe}{\end{defeni}}

\newenvironment{remark}[1]%
    { \begin{Rem} \label{R:#1}  \ifShowLabels \TeXref{R:#1} \fi }%
    { \end{Rem} }
\newcommand{\rem}[1]{\begin{remark}{#1}}
\newcommand{\erem}{\end{remark}}

\newenvironment{conjec}[1]%
    { \begin{Conj} \label{Co:#1}  \ifShowLabels \TeXref{Co:#1} \fi }%
    { \end{Conj} }
\renewcommand{\conj}[1]{\begin{conjec}{#1} \sl }
\newcommand{\econj}{\end{conjec}}

\newcommand{\eq}[1]%
    { \ifShowLabels \TeXref{E:#1} \fi
       \begin{equation} \label{E:#1} }
\newcommand{\eeq}{ \end{equation} }

\newcommand{\prf}{ \begin{proof} }
\newcommand{\epr}{ \end{proof} }

\newcommand\alp{\alpha}		
		
\newcommand\gam{\gamma}		
\newcommand\del{\delta}		\newcommand\Del{\Delta}
\newcommand\eps{\varepsilon}

\newcommand\kap{\kappa}
\newcommand\lam{\lambda}		\newcommand\Lam{\Lambda}
\newcommand\sig{\sigma}		

\newcommand\ome{\omega}		

\newcommand\calA{{\mathcal{A}}}
\newcommand\calB{{\mathcal{B}}}
\newcommand\calC{{\mathcal{C}}}

\newcommand\calE{{\mathcal{E}}}

\newcommand\calH{{\mathcal{H}}}

\newcommand\calJ{{\mathcal{J}}}
\newcommand\calK{{\mathcal{K}}}

\newcommand\calM{{\mathcal{M}}}

\newcommand\calO{{\mathcal{O}}}

\newcommand\calS{{\mathcal{S}}}

\newcommand\calZ{{\mathcal{Z}}}

		\newcommand\bfB{{\mathbf B}}

		\newcommand\bfG{{\mathbf G}}

		\newcommand\bfM{{\mathbf M}}

		\newcommand\bfP{{\mathbf P}}

		\newcommand\bfT{{\mathbf T}}
		\newcommand\bfU{{\mathbf U}}

		\newcommand\bfX{{\mathbf X}}
		
		\newcommand\bfZ{{\mathbf Z}}

\newcommand\RR{\mathbb{R}}

\renewcommand\AA{\mathbb{A}}

\newcommand\ZZ{\mathbb{Z}}

\newcommand\CC{\mathbb{C}}

	\newcommand\grg{{\mathfrak{g}}}

	\newcommand\grm{{\mathfrak{m}}}

	\newcommand\grp{{\mathfrak{p}}}

	\newcommand\gru{{\mathfrak{u}}}

\newcommand\sdp{\times \hskip -0.3em {\raise 0.3ex
\hbox{$\scriptscriptstyle |$}}} 

\newcommand\End{\operatorname{End\,}}

\newcommand\Hom{\operatorname {Hom}}

\newcommand\Ind{\operatorname{Ind}}

\newcommand\oX{{\overline{X}}}

\newcommand\tilG{{\widetilde{G}}}

\newcommand\tilT{{\widetilde{T}}}

\newcommand\tilX{{\widetilde{X}}}

\newcommand\x{\times}
\newcommand\ten{\otimes}

\newcommand{\ra}{\rangle}
\newcommand{\la}{\langle}

\newcommand\tK{\textsf{K}}
\newcommand\val{\operatorname{val}}
\newcommand\aff{\operatorname{aff}}
\newcommand\vol{\operatorname{vol}}
\begin{document}

\title{Remarks on the asymptotic Hecke algebra}
\author{Alexander Braverman and David Kazhdan}
\begin{abstract}Let $G$ be a split reductive $p$-adic group,
$I\subset G$ be an Iwahori subgroup, $\calH(G)$ be the Hecke algebra and $\calC(G)\supset \calH(G)$ be
the Harish-Chandra Schwartz algebra. The purpose of this note is to define (in spectral terms)  a subalgebra $\calJ(G)$ of $\calC(G)$, containing $\calH(G)$, which we consider as an algebraic version of $\calC(G)$.
We show that the subalgebra $\calJ(G)^{I\times I}\subset \calJ(G)$ is isomorphic to the Lusztig's asymptotic Hecke algebra $J$
 and explain a relation between the algebra $\calJ(G)$ and the Schwartz space of the basic affine space studied in \cite{BK}.
\end{abstract}
\address{A.B.: Department of Mathematics, University of Toronto, Perimeter Institute for Theoretical Physics and Skolkovo Institute for Science and Technology;}
\address{D.K.: Department of Mathematics, Hebrew University, Jerusalem, Israel}

\maketitle

\sec{}{Introduction and statement of the results}
\ssec{}{Notation}Let $F$ be a non-archimedian local field with ring of integers $\calO$; we shall choose a generator $\pi$ of the maximal ideal of $\calO$. Typically, we shall denote algebraic varieties over $F$ with boldface letters (e.g. $\bfG,\bfX$ etc.) and the corresponding sets of $F$-points -- with corresponding Roman letters (e.g. $G,X$ etc.).

In what follows we fix a connected split reductive group $\bfG$ over $F$ with a Borel subgroup $\bfB$ it unipotent radical $\bfU$, maximal split torus $\bfT=\bfB/\bfU$. Let $\Lam$
be the lattice of cocharacters of $\bfT$
 and $\Lam ^\vee$ be the lattice of characters of $\bfT$.

 We write $K_0=G(\calO)$ and denote by $I\subset K_0$ an Iwahori subgroup of $G$. We denote by $\calH(G)$ the full Hecke algebra of $G$ and by $\calH(G,I)$ the Iwahori-Hecke subalgebra; we shall also denote by $\calH_{\aff}$ the corresponding algebra over $\CC[v,v^{-1}]$ (thus $\calH(G,I)$ is obtained from $\calH_{\aff}$ by specializing $v$ to $q^{1/2}$).
We denote by $\calM(G)$ the category of smooth representations of $G$. This is the same as the category of locally unital modules over $\calH(G)$. For any  smooth representation $(\pi ,V)$ we denote by  $(\pi ^\vee  ,V^\vee )$ the subrepresentation of smooth vectors in the representation of $G$ on the space of linear functionals on $V$. For any $v\in V,\lam \in V^\vee$ we denote by $m_{v, \lam}:G\to \CC$ the matrix coefficient
 $m_{v, \lam}(g):=\lam (\pi (g)(v)).$

\ssec{}{Matrix Paley-Wiener theorem}Let $P$  be a parabolic subgroup  of $G$ with a Levi group $M$. The set  $X_M$ of unramified characters of $M$ is
equal to $\Lam ^\vee _M\otimes \CC ^\times$ where
 $\Lam ^\vee _M\subset \Lam^{\vee}$ is the subgroup of characters of $\bfT$ trivial on $\bfT \cap [\bfM,\bfM]$. So
 $X_M$ has a structure of a complex algebraic variety; the algebra of polynomial functions on $X_M$ is equal to $\CC[\Lam_M]$ where $\Lam_M$ is the lattice dual to $\Lam_M^{\vee}$. We denote by $X_{M,t}\subset X_M$ the subset of unitary characters.

 For any $(\sig ,V)\in\calM(M)$ we denote by $i_{GP}(\sig)$ the corresponding unitarily induced object of $\calM(G)$. As a representation of $K_0$ this representation is equal to $\text{ind} ^{K_0}_{P\cap K_0}(\sig)$. So for any  unramified character $\chi:M\to \CC^*$ the space of the representation
$i_{GP}(\sig\ten \chi)$  is isomorphic to the space $V_\chi$  of the representation $i_{GP}(\sig )$ and  is independent on a choice of $\chi$.
Since  $X_M$ has a structure of an algebraic variety over $\CC$ it make sense to say that a family $\eta_\chi \in \text{End} (V_\chi ), \chi \in \bfT ^\vee$ is a regular (or a smooth) function of $\chi$.

We denote by $Forg :\calM (G)\to Vect$ the forgetful functor,
 by $\widetilde{\calE(G)}=\{ e(\pi )\}$ the ring of endomorphisms of $Forg$ and define  $\calE(G)\subset \widetilde{\calE(G)}$ as the subring of endomorphisms $\eta_{\pi}$ such that

1) For any
 Levi subgroup $M$ of $G$ and $\sig \in Ob( \calM (M))$,
 the endomorphisms $\eta_{i_{GP}(\sig\ten \chi)}$  are regular functions of $\chi$.

2) There exists an open compact subgroup $K$ of $G$ such that $\eta_{\pi}$ is $K\times K$-invariant for every $\pi$.

By definition, we have a homomorphism
$$
PW :\calH (G)\to \calE(G),\quad f\mapsto \pi (f).
$$
The following is usually called "the matrix Paley-Wiener theorem" (cf. \cite{Ber}, Theorem 25):
\th{pw}The map $PW$ is an isomorphism.
\eth

The group $G\times G$ acts on $\calE(G)$ in the obvious way. We denote by $\calE^I(G)\subset \calE(G)$  the subring $I\times I$-invariant elements.
It is clear that we can interpret $\calE^I(G)$ as a subring of the ring of endomorphisms of the forgetful functor $Rep( \calH (G,I))\to Vect$.

\ssec{}{Harish-Chandra algebra}
Recall that for any $g\in G$, there exists a unique dominant coweight $\lam(g)$ of $\bfT$ such that $g\in \bfG(\calO)\pi^{\lam(g)}\bfG(\calO)$.
Let us set $\Del(g)=q^{\la\lam,\rho\ra}$. Then we say that a function $f:G\to \CC$ is a Schwartz function if

\medskip
(a) There exists an open compact subgroup $K$ of $G$ such that $f$ is both left and right $K$-invariant.

(b) For any polynomial function $p:G\to F$ and $n>0$,  there exists a constant $C=C_{p,n}\in \RR_{>0}$ such that
$$\Del(g)|f(g)|\leq C\ln ^{-n}(1+| p(g)|)$$
for all $g\in G$.

\medskip
\noindent
We denote by $\calC(G)$ the space of all Schwartz functions.
It is known that $\calC (G)$ has  an algebra structure  with respect to convolution (cf. \cite{sil}, Chapter 4 or \cite{Wal}).

Obviously we have
$\calH(G)\subset \calC(G)$.

For an open compact subgroup $K$ of $G$ we denote by $\calC(G,K)$ the space of $K\x K$-invariants in $\calC(G)$.

 Below we recall the spectral description of $\calC(G)$.
\ssec{}{Tempered representations}Let $(\pi ,V)$ be a representation of $G$  of finite length with central character $\eta:Z(G)\to \CC^*$. Recall that $\pi$ is called tempered if

\medskip
(1) $\pi$ is unitary. In particular, $\eta$ is unitary (i.e., it takes values in $S^1\subset \CC^*$).
In this case the absolute values $|m_{v,\lam}|$ of  matrix coefficents of $V$ are functions on $G/Z(G)$.

(2) For any $\eps>0$  and any matrix coefficient $m_{v,\lam}$ of $\pi$ we have

 $$|m_{v,\lam}|\in L^{2+\eps}(G/Z(G)).$$

\medskip
\noindent

The following facts are well-known (cf. \cite{sil}, \cite{Wal}):

(F1) Let $\pi$ be a tempered representation of $G$. Then the action of $\calH(G)$ extends naturally to an action of $\calC(G)$.

(F2) Let $P$ be a parabolic subgroup of $G$ with a Levi group $M$. Let $\sig$ be a tempered irreducible representation of $M$. Then the representation $i_{GP}(\sig)$ is tempered.

(F3) For a generic unitary character $\chi:M\to S^1$ the representation $i_{GP}(\sig\ten \chi)$ (which is tempered by F2) is irreducible.

\noindent

We denote by $\calM _t(G)\subset \calM (G)$ the subcategory of  tempered representations.
As follows from (F2), for any tempered representation $\sig$ of $M$ and a unitary character $\chi$ of $M$
 the representations  $i_{GP}(\sig\ten \chi)$ belong to $\calM _t(G)$.

 Let $\calE_t(G)$ be the subring of endomorphisms $\{ \eta\}$ of the forgetful functor $Forg_t:\calM _t(G)\to Vect$ such that

(1t) $\eta_{i_{GP}(\sig\ten \chi)}$ is a smooth function of $\chi \in X_{M,t}$
for any Levi subgroup $M$ of $G$ and $\sig \in Ob( \calM _t(M))$.

(2t) $\eta$ is $K\times K$-invariant for some open compact subgroup of $G$.

\noindent
The fact F1 can be upgraded to the following version of the matrix Paley-Wiener theorem (cf. \cite{Wal}):

\th{pw-tempered}
The map $f\mapsto \pi(f)$ defines an isomorphism between $\calC(G)$ and $\calE_t(G)$.
\eth

We denote by $\calE_t^I\subset \calE_t(G)$  the subring of $I\times I$-invariant elements of $\calE_t(G)$.

\ssec{}{Asymptotic Hecke algebra}

Recall that we denote by $\calH_{\aff}$ the "algebraic" version of $\calH(G,I)$ which is an algebra over $\CC[v,v^{-1}]$. Let us assume that $\bfG$ is of adjoint type. In \cite{Cells2}, G.~Lusztig defined the so-called
asymptotic Hecke algebra $J$ (we are going to recall the definition in \refs{proof}). This is an algebra over $\CC$\ \footnote{In fact $J$ can be defined over $\ZZ$ but to simplify the notation we shall always work over $\CC$} and there is a canonical embedding $\calH(G,I)\hookrightarrow J\ten \CC[v,v^{-1}]$ which becomes an isomorphism after some completion. Moreover, one can show that the specialization of this embedding to any $q\in \CC^*$ is also injective. Hence we get an embedding
$\calH(G,I)\hookrightarrow J$.

One of the main purposes of this paper is to formulate and prove a version of matrix Paley-Wiener theorem for $J$.

Let $P$ be a parabolic subgroup with Levi group $M$. We say that an unramified character $\chi:M\to \CC^*$ is (non-strictly) positive if for any coroot $\alp$ of $\bfG$, such that the corresponding root subgroup lies in the unipotent radical $\bfU_{\bfP}$ of $\bfP$ (which in particular defines a homomorphism $\alp:F^*\to Z(M)$), we have $|\chi(\alp(x))|\geq 1$ for $|x|\geq 1$.

Let $\calE_{\calJ}(G)$ be ring of collections $\{ \eta_{\pi}\in \End_{\CC}(V)|\ \text{for tempered irreducible $(\pi,V)$}\}$ which extend to a rational function $E_{i_{GP}(\sig\ten \chi)}\in \End_{\CC}(\sig\ten \chi)$ for every tempered irreducible representation $\sig$ of $M$  and which are

a) regular on the set  of  characters $\chi$
such that $\chi^{-1}$ is (non-strictly) positive.

b) $K$-invariant for some open compact subgroup $K$ of $G$.

\noindent
As follows from the definition, we have an embedding $\calE_{\calJ}^I(G)\to \calE^I_t(G)$.
\th{main}
\begin{enumerate}
\item
Let $(\pi,V)$ be an irreducible tempered representation of $G$.
Then the action the action of $\calH(G,I)$ on $V^I$ extends uniquely to $J$.
\item
Let $P$ be a parabolic subgroup of $G$ with Levi group $M$, let $\sig$ be an irreducible tempered representation of $M$, and let $\chi$ be a (non-strictly) positive character of $M$ and $(\pi,V)=i_{GP}(\sig\ten \chi^{-1})$. Then the action of $\calH(G,I)$ on $V^I$ extends to an action of $J$. This extension is unique up to isomorphism.
 \item
 The map $f\mapsto \pi(f)$ defines an isomorphism between $J$ and $\calE_{\calJ}^I(G)$.

 \end{enumerate}
\eth
It follows immediately from \reft{pw}, \reft{pw-tempered} and \reft{main} that we have inclusions $\calH(G,I)\subset J\subset \calC(G,I)$.
\reft{main} allows giving the following
\defe{} We define $\calJ(G)$ to be the preimage of $\calE_{\calJ}(G)$ in $\calC(G)$. Note that we have natural embeddings $\calH(G)\subset \calJ(G)\subset \calC(G)$.
\edefe

The algebra $\calJ(G)$ can be thought of as a "beyond Iwahori" version of $\calJ(G,I)$. It follows again from \reft{pw} and \reft{pw-tempered} that we have the embeddings $\calH(G)\subset \calJ(G)\subset \calC(G)$. Also, this definition makes sense for any reductive $G$.

The proof of \reft{main} is given in the next Section; it is essentially an exercise on manipulating the results from \cite{Cells2}, \cite{Cells3} and \cite{Cells4}.
In \refs{schwartz} we also explain a connection between the algebra $\calJ(G)$ and the Schwartz space of the basic affine space studied in  \cite{BK}. In the appendix we give some examples of elements of $J$ (viewed as functions on $G$) for $G=SL(2,F)$.

\ssec{}{An algebraic version}
Let us explain a version of \reft{main} which is algebraic in the sense that it doesn't use the notion of "positive real number". To emphasize this we are going to work not over $\CC$ but over an arbitrary algebraically closed field $\texttt{K}$ of characteristic $0$.
First, let us make the following definitions:

(1) A character $\chi:F^*\to \texttt{K}^*$ is called special if $\chi(a)=q^{\frac{1}{2}r\cdot  \text{val}(a)}$ with $r\in \ZZ_{>0}$

(2) A representation of $F^*$ on a  $\texttt{K}$-vector space $V$ is called special if there is a non-zero vector in $V$ on which $F^*$ acts by a special character.

(3) Let $G,P,M$ be as above. Then a representation $\sig$ of $M$ is called special if there exists a positive coroot $\alp:F^*\to Z(M)$ whose composition with $\sig$ is special

(4) A representation $\pi$ of $G$ is called quasi-tempered if for any parabolic $P$ with Levi decomposition $P=MU$, the Jacquet functor $r_{GP}(\pi)$ is not a special representation of $M$.

(5) Let $\calE_{\calJ,\texttt{K}}$ be the algebra
$\{ \eta_{\pi}\in \End_{\CC}(V)|\ \text{for quasi-tempered irreducible $(\pi,V)$}\}$ which extend to a rational function $E_{i_{GP}(\sig)}\in \End_{\CC}(\sig)$ for every non-special irreducible representation $\sig$ of $M$. For an open compact subgroup $K$ of $G$ we denote by $\calE{\calJ,\texttt{K}}(G,K)$ the algebra of $K\times K$-invariant elements in $\calE_{\calJ,\texttt{K}}$.

\medskip
\noindent
With these definitions in mind we have
\th{main'}
Let $J_{\texttt{K}}=J_{\ZZ}\ten \texttt{K}$.
Then
\begin{enumerate}
\item
Let $(\pi,V)$ be an irreducible quasi-tempered representation of $G$. Then the action of $\calH(G,I)$ on $V$ extends uniquely to $J_{\texttt{K}}$.
\item
Let $P$ be a parabolic subgroup of $G$ with Levi group $M$. Let $\sig$ be a non-special representation of $M$. Let $(\pi,V)=i_{GP}(\sig)$. Then the action of $\calH(G,I)$ on $V^I$ extends uniquely to $J_{\texttt{K}}$.
\item
The map $f\mapsto \pi(f)$ defines an isomorphism between $J_{\texttt{K}}$ and $\calE_{\calJ,\texttt{K}}(G,I).$
\end{enumerate}
\eth
It is not difficult to show that  $\calE_{\calJ,\texttt{K}}$ is naturally isomorphic to some algebra of locally constant functions on $G$ (where the algebra structure is given by convolution) which contains $\calH_{\texttt{K}}(G)$; we shall denote this algebra by $\calJ_{\texttt{K}}(G)$.
We shall not pursue the details in this paper.

\ssec{}{Acknowledgements}We would like to thank R.~Bezrukavnikov, G.~Lusztig and V.~Ostrik for numerous discussions and help during the preparation of this paper.

The project has received funding from ERC under grant agreement 669655. In addition the first-named author was partially supported by NSERC.
\sec{proof}{Proof of \reft{main}}
In this Section we shall assume that $\bfG$ is semi-simple and adjoint. Then the affine Weyl group $W_{\aff}=\Lam\rtimes W$ is in general not a  Coxeter group, but it still has a length function $\ell$ defined e.g. in Section 1.1 of \cite{Cells4}.
\ssec{}{Definition of $J$} We are keeping the notations of \cite{Cells2} and \cite{Cells4}.
Let $\calA=\CC[v,v^{-1}]$ and $\calA^+=\CC[v]$. We let as before $\calH_{\aff}$ denote the corresponding affine Hecke algebra $\calH_{\aff}$. It has bases $\{ T_w,\tilT_w,C_w\}_{w\in W_{\aff}}$ where $T_w$ is the standard basis of $\calH_{\aff}$, $\tilT_w=v^{-\ell(w)}T_w$ and $C_w$ is the Kazhdan-Lusztig basis which satisfies
$$
C_w=\sum\limits_{y\leq w} (-1)^{\ell(w)-\ell(y)} v^{\ell(w)-\ell(y)} P_{y,w}(v^{-2})\tilT_y,
$$
where $P_{y,w}$ is a polynomial.

For $x,y\in W_{\aff}$, let us set
$$
C_xC_y=\sum\limits_{z\in W_{\aff}} h_{x,y,z} C_z.
$$

For any $z\in W_{\aff}$, we let $a(z)$ be the smallest integer such that $u^{a(z)}h_{x,y,z}\in \calA^+$ for all $x,y$. Let $\gam_{x,y,z}\in \ZZ$ be such that $u^{a(z)}h_{x,y,z}-\gam_{x,y,z}\in u\calA^+$. Then let $J$ denote the $\CC$-vector space with basis $t_w, w\in W_{\aff}$; it has a ring structure defined by
$$
t_xt_y=\sum\limits_{z\in W_{\aff}} \gam_{x,y,z}t_z.
$$

Let $k=\CC((v^{-1}))$ with the natural $v^{-1}$-adic topology on it. Let also $\widehat{\calH_{\aff}}$ be the completion of $\calH_{\aff}$ consisting of all (possibly infinite) sums
$\sum b_x C_x$ with $b_x\in \widehat{\calA}$ such that $b_x\to 0$ when $\ell(x)\to \infty$ (it is clear that one can replace $C_x$ by $\tilT_x$ in the above definition).
In Section 2.4 of \cite{Cells2}, G.~Lusztig defines a homomorphism $\calH_{\aff}\hookrightarrow J\ten \calA$ which gives rise to an isomorphism
between $\widehat{\calH_{\aff}}$ and $J\widehat{\ten} k$ (the latter denotes some completed tensor product). Also, the above embedding makes $J\ten \calA$ into a finitely generated module over $\calH_{\aff}$.
\ssec{reph}{Representations of $\calH_{\aff}$ and $\calH(G,I)$} Set $J_{\calA}=J\ten \calA, J_{k}=J\ten k, \calH_{\aff,k}=\calH_{\aff}\underset{\calA}\ten k$. It follows from the above that $\calH_{\aff,k}$ can be regarded as a subalgebra of $J_{k}$.
In what follows we denote by $G^{\vee}$ the Langlands dual group of $\bfG$ (over $\CC$).

Let $(s,u,\rho)$ be a triple where

(a) $s\in G^{\vee}(\CC)$ is a semi-simple element,

(b) $u\in G^{\vee}(\CC)$ is a unipotent element such that $su=us$,

(c) $\rho$ is an irreducible representation of the group of components of the centralizer $Z_{G^{\vee}}(s,u)$ of the pair $(s,u)$.

Recall that in \cite{KL} Kazhdan and Lusztig define a representation $\tK(s,u,\rho)$ of $\calH_{\aff}$. This representation may be 0; we shall say that the triple $(s,u,\rho)$ is admissible if $\tK(s,u,\rho)\neq 0$. We shall denote $\tK(s,u,\rho,q)$ the specialization of $K(s,u,\rho)$ to $v=q^{1/2}$ viewed as a representation of $\calH(G,I)$. Moreover, if $\rho$ is some (not necessarily irreducible) representation of the group of components of $Z_{G^{\vee}}(s,u)$, we denote by $\tK(s,u,\rho,q)$ the direct sum of $\tK(s,u,\rho_i)$ where $\rho$ is the direct sum of irreducible representations $\rho_i$).
Then the following facts are true:
\begin{enumerate}
\item (cf. Theorem 7.12 in \cite{KL})
  $\tK(s,u,\rho,q)$ has a unique simple quotient which we shall denote by $L(u,s,\rho,q)$.
  \item (cf. Theorem 8.2 in \cite{KL})
  Assume that $s$ is compact. Then the corresponding representation $\tK(s,u,\rho,q)$ of $\calH(G,I)$ is tempered and irreducible. Moreover, every irreducible tempered representation of $\calH(G,I)$ is isomorphic to $\tK(s,u,\rho,q)$ for a unique admissible triple $(s,u,\rho)$ with compact $s$.
  \item (cf. Theorem 6.2 in \cite{KL})
  Let $P$ be a parabolic subgroup of $G$ with Levi group $M$ and let $P^{\vee}$ and $M^{\vee}$ be the corresponding parabolic and Levi subgroups in $G^{\vee}$. Let also $\grp^{\vee}, \grm^{\vee}$ be their Lie algebras. Assume that $s,u\in M^{\vee}$. Let $Z_{M^{\vee}}(s,u)$ be the centralizer of $(s,u)$ in $M^{\vee}$ and let $Z_{G^{\vee}}(s,u)$ be the centralizer of $(s,u)$ in $G^{\vee}$.

   Now let $(\grg^{\vee}/\grp^{\vee})_u$ be the kernel of $1-u$ on $\grg^{\vee}/\grp^{\vee}$. Assume that $q^{-i}$ is not an eigen-value of $s$ on $(\grg^{\vee}/\grp^{\vee})_u$ for every $i>0$. Let $\tK_M(s,u,\rho_j,q)$ denote the corresponding representation of $M$. Then $i_{GP}(\tK_M(s,u,\rho,q))$ is isomorphic to $\tK(s,u,\widetilde{\rho},q)$ where $\widetilde{\rho}=\Ind_{Z_M(s,u)/Z^0_M(s,u)}^{Z_G(s,u)/Z^0_G(s,u)}\rho$. In particular, this is true if $s$ is of the form $s'\cdot \chi^{-1}$ where $s'$ is compact and $\chi$ is a (non-strictly) positive character of $M$.
  \end{enumerate}

In \cite{Cells4} Lusztig proves the following result:
\th{lus-main}Let $E$ be an irreducible representation of $J$ and let $E_{k}=\rho\underset{\CC}\ten k$.
Then there exists a unique triple $(s,u,\rho)$ such that $E_{k}|_{\calH_{\aff,k}}$ is isomorphic to $\tK(s,u,\rho)\underset{\CC}\ten k$ (we shall denote the latter $\calH_{\aff,k}$-module by $\tK(s,u,\rho)_k$). Moreover, every admissible triple $(s,u,\rho)$ arises in this way.
\eth
\reft{lus-main} implies that we have a bijection between irreducible representations of $J$ and admissible triples $(s,u,\rho)$. For any such triple $(s,u,\rho)$ we shall denote by $E(s,u,\rho)$ the corresponding irreducible representation of $J$.

Note that by specializing the embedding $\calH_{\aff}\hookrightarrow J\ten \calA$ to $v=q^{1/2}$ we get a homomorphism $\calH(G,I)\to J$ which is injective by Proposition 1.7 of \cite{Cells3}.
We now claim the following:
\th{foremb}
\begin{enumerate}
\item
Let $\pi$ be an irreducible tempered representation of $\calH(G,I)$. Then $\pi$ extends uniquely to an irreducible representation of $J$.
\item
Any module  of the form $i_{GP}(\sig\ten \chi^{-1})$, where $\sig$ is an irreducible tempered representation of the Levi group $M$ and $\chi$ is a (non-strictly) positive character of $M$, extends to an irreducible module over $J$. This extension is unique up to isomorphism.
\item
Let $M$ be a Levi subgroup of $G$ and let $(s,u,\rho)$ be an admissible triple for in $M^{\vee}$ with compact $s$. Then there exists
a $J\ten \CC[\Lam_M]$-module $\calM(s,u,\rho)$ whose fiber at any non-strictly positive $\chi$ is isomorphic to the $J$-module from (2).
\end{enumerate}
\eth

In order to prove \reft{foremb} we shall use the following result of N.~Xi \cite{xi}:
\th{xi}Let $(s,u,\rho)$ be an admissible triple. Then the restriction of $E(s,u,\rho)$ to $\calH(G,I)$ has a unique irreducible quotient, which is isomorphic to $L(s,u,\rho,q)$. Moreover, any irreducible subquotient of
the kernel of the map $E(s,u,\rho)|_{\calH(G,I)}\to L(s,u,\rho,q)$ is not isomorphic to $L(s,u,\rho,q)$.
\eth

This result implies the following:
\cor{corxi}
The representation $E(s,u,\rho)|_{\calH(G,I)}$ of $\calH(G,I)$ is isomorphic to $\tK(s,u,\rho,q)$.
\ecor
\prf
By definition we have an isomorphism $E(s,u,\rho)|_{\calH(G,I)}\simeq \tK(s,u,\rho)_k$. This obviously implies that

(a) $\dim E=\dim_k \tK(s,u,\rho)_k=\dim _{\CC}\tK(s,u,\rho,q)$.

(b) There exists a non-zero homomorphism $\tK(s,u,\rho,q)\to E(s,u,\rho)|_{\calH(G,I)}$.

\noindent
According to \cite{KL}, Theorem 7.12, the representation $\tK(s,u,\rho,q)$ has a unique simple quotient $L(s,u,\rho,q)$.
It now follows from \reft{xi} that the map $\tK(s,u,\rho,q)\to E(s,u,\rho)|_{\calH(G,I)}$ is surjective (indeed, otherwise it would land inside the maximal proper submodule of $E(s,u,\rho)|_{\calH(G,I)}$ which doesn't contain $L(s,u,\rho,q)$ as a subquotient. On the other hand, any non-zero image of $\tK(s,u,\rho,q)$ has a quotient isomorphic to $L(s,u,\rho,q)$). Since the dimensions of these two modules are equal, it follows that
the map $\tK(s,u,\rho,q)\to E(s,u,\rho)|_{\calH(G,I)}$ is an isomorphism.
\epr
\ssec{}{Proof of \reft{foremb}}
Let $\pi$ be an irreducible tempered representation of $\calH(G,I)$. Then $\pi$ is isomorphic to $\tK(s,u,\rho,q)$ for some admissible triple $(s,u,\rho)$ with compact $s$. Now \refc{corxi} implies that $\pi$ extends to an irreducible representation of $J$ and any two such extensions
are isomorphic as abstract $J$-modules. In other words, any two such extensions are conjugate by means of some $\text{Aut}_{\calH(G,I)}(\pi)=\CC^{\times}$. Hence any two such extensions are equal. This proves assertion (1).

Assertion (2) follows in a similar way from statement (3) above \reft{lus-main}.
Finally, the module $\calM(s,u,\rho)$ is constructed in the following way.
Let $\calB^{\vee}$ denote the flag variety of $G^{\vee}$ and let $\calB^{\vee}_{s,u}$ denote the variety of $(s,u)$-fixed points on $\calB^{\vee}$. Then it follows from \cite{Cells4} that
$J$ acts on the equivariant $K$-theory $K_{Z(M^{\vee})}(\calB^{\vee}_{s,u})$ (here $Z(M^{\vee})$ denotes the center of $M^{\vee}$).
Moreover, this $K$-theory has a natural action of the centralizer $Z_{M^{\vee}}(s,u)$ and we let $\calM(s,u,\rho)$ denote its $\rho$-isotypic component with respect to $Z_{M^{\vee}}(s,u)$. It easily follows from the above that $\calM(s,u,\rho)$ satisfies the requirements of (3).

\ssec{}{Proof of \reft{main}}
The first two assertions of \reft{main} are exactly the two assertions of \reft{foremb}.
So, it remains to prove the 3rd assertion.

It follows from part (3) of \reft{foremb} that the map $J\to \calE^I_{\calJ}(G)$ is well defined. We now need to prove that it is an isomorphism.
First, given $h\in J$ for any $P,M,\sig$ as above, we can define $\eta_h(\chi)\in \End_{\CC}(i_{GP}(\sig\ten \chi^{-1}))$ where $\chi$ is a positive unramified character of $M$;
the fact that $\eta_h(\chi)$ depends rationally on $\chi$ follows immediately from the fact that the embedding $\calH(G,I)\hookrightarrow J$ is an isomorphism over the generic point of the center of $\calH(G,I)$.
Thus we get an injective map $J\hookrightarrow \calE_{\calJ}^I(G)$. We now want to prove that this map is also surjective.

For a unipotent element $u$ in $G^{\vee}$ let $\calZ_u$ denote the algebra of ad-invariant polynomial functions on the centralizer $Z_{G^{\vee}}(u)$ of $u$ in $G^{\vee}$. Let $\calZ$ denote the direct sum of all the $\calZ_u$ (where the sum is taken over
conjugacy classes of unipotent elements in $G^{\vee}$). Then $\calZ$ maps to the center of both
$J$ and $\calE^I_{\calJ}(G)$ and both algebras are finitely generated modules over $\calZ$.

To prove that the desired surjectivity holds, it is enough to prove that it holds modulo every maximal ideal of $\calZ$. Let $\grm$ be such an ideal. Set $J_{\grm}=J/\grm J, \calE_{\grm}=\calE_{\calJ}^I(G)/\grm \calE_{\calJ}^I(G)$. It is enough to prove that the map $J_{\grm}\to \calE_{\grm}$ is surjective for every $\grm$.
Let $(\pi_1,V_1),\cdots, (\pi_n,V_n)$ be all the different (non-isomorphic) representations of $G$ which have the form $i_{GP}(\sig\ten \chi^{-1})$ with tempered $\sig$ and positive $\chi$ such that $V_i^I\neq 0$ and  $\calZ$ acts on $V_i^I$ through  the quotient by $\grm$. Then by definition $\calE_{\grm}$ embeds into $\oplus_i \End_{\CC}(V_i^I)$. On the other hand, $V_1^I,\cdots, V_n^I$ are non-isomorphic irreducible representations of
$J_{\grm}$ and hence the map $J_{\grm} \to \oplus_i \End_{\CC}(V_i^I)$ is surjective.

\sec{schwartz}{Connection to the Schwartz space of $G/U$}

\ssec{dig}{Digression on \cite{BK}}
Let $\bfU$ be a maximal unipotent subgroup of $\bfG$ defined over $F$, $U=\bfU(F)$. Set $X=G/U$; it is endowed with a natural action of $G\times T$.
Let us denote by $\calS_c(X)$ the space of
locally constant compactly supported functions on $X$, also let $C^{\infty}(X)$ denote just the space functions $f: X\to \CC$, such that there exists an open compact subgroup $K$ of $G$ such that $f$ is $K$-invariant.
Let $G\times T$ acts on these spaces is such a way that the action of $G$ comes from the right action of $G$ on $X$ and the action of $T$ comes from the right action of $T$ on $X$ twisted by the character
$t\mapsto q^{\la \val(t),\rho\ra}$ of $T$ where $\val:T\to \Lam$ denotes the natural homomorphism. In \cite{BK} we have defined the  Schwartz space $\calS(G/U)$ of functions on the basic affine space $G/U$ which contains $\calS_c(G/U)$ and it is contained in $C^{\infty}(X)$ in the case when $\bfG$ is simply connected. Let us recall this definition.

The space $X$ has unique up to scalar $G$-invariant measure and we denote by $L^2(X)$ the $L^2$-space with respect to this measure.
When $\bfG$ is simply connected one can construct a natural action of the Weyl group $W$ on $L^2(X)$ by unitary operators $\Phi_w$ which commute with $G\times T$.
In order to define these operators it is enough to consider the case when $w=s_{\alp}$ -- a simple reflection (here $\alp$ is a simple root of $G$).
Let us recall this definition as it will be used in the future.

For a simple root $\alpha$ let $\bfP_\alpha\subset \bfG$ be
the minimal parabolic
of type $\alpha$ containing $B$. Let $\bfB_\alpha$ be the commutator
subgroup of $P_\alpha$,
and denote $\bfX_\alpha:=\bfG/\bfB_\alpha$. We have an obvious projection
of homogeneous spaces $\pi _\alpha:
\bfX\rightarrow \bfX_\alpha$.
It is a fibration with the fiber $\bfB_\alpha /\bfU=\AA^2-\{0\}$.

Let $\overline \pi _\alpha :{\overline \bfX^\alpha} \to \bfX_\alpha$
be the relative
affine completion of the morphism $\pi _\alpha$. (So  $\overline \pi _\alpha$
is the affine morphism corresponding to the sheaf of algebras $\pi_{\alpha *}
(\calO_{\bfX})$ on $\bfX_\alpha$.) Then  $\overline \pi _\alpha$
has the structure
of a 2-dimensional vector bundle; $\bfX$ is identified with the complement to
the zero-section in $\overline \bfX^\alpha$. The $\bfG$-action on $\bfX$
obviously
extends to $\overline \bfX^\alpha$; moreover, it is easy to see that the
determinant of the vector bundle $\overline \pi _\alpha$ admits a canonical
(up to a constant) $\bfG$-invariant trivialization, i.e.
$\overline \pi _\alpha$ admits unique up to a constant $\bfG$-invariant
fiberwise symplectic form $\ome_\alp$.
We will fix such a form for every $\alp$.

Obviously $L^2(X)=L^2({\overline X^\alpha})$. Thus we define
$\Phi_\alp=\Phi_{s_\alp}$ to be equal to the Fourier transform in the fibers
of $\overline \pi _\alpha$, corresponding to the identification
of ${\overline X^\alpha}$ with the dual bundle
by means of $\ome_\alp$.

Then
$$
\calS(X)=\sum\limits_{w\in W} \Phi_w(\calS_c(X)).
$$
We can extend the above definition to the case when $\bfG$ is not necessarily simply connected. First the definition of \cite{BK} works without any change in the case when $[\bfG,\bfG]$ is simply connected.
Now, given any connected reductive $\bfG$ there always exists an algebraic reductive group $\widetilde{\bfG}$ and a central torus $\bfZ$ in $\widetilde{\bfG}$ so that $\bfG=\widetilde{\bfG}/\bfZ$.
We now denote by $\tilX$ the basic affine space for $\tilG=\widetilde{\bfG}(F)$ and we set $\calS(X)=\calS(\tilX)^T$.
With this definition most results of \cite{BK} extend word-by-word to any $\bfG$.

\ssec{}{Action of $\calJ(G)$ on $L^2(X)$}
By definition we have $\calS(X)\subset L^2(X)$. We claim that $\calC(G)$ acts on $L^2(X)$.
Indeed, we have
$$
L^2(X)=\bigoplus\limits_{\theta:\bfT(\calO)\to S^1} L^2(X)_{\theta},
$$
where $L^2(X)_{\theta}$ denotes the subspace of $L^2(X)$ on which $\bfT(\calO)$ acts by $\theta$. Now, each $L^2(X)_{\theta}$ is a direct integral of $G$-representations of the form $i_{GB}(\chi)$
where $\chi$ is a unitary character of $T$ over a compact base (isomorphic to $(S^1)^{\dim \bfT}$) and therefore it acquires a natural action of $\calC(G)$ (since it acts on each $i_{GB}(\chi)$ with unitary $\chi$).

In particular, the algebra $\calJ(G)$ acts on $L^2(X)$. The following conjecture provides an alternative definition of $\calS(X)$.

\conj{main-sch}
We have $\calS(X)=\calJ(G)\cdot \calS_c(X)$.
\econj

\noindent
{\bf Remark.} We claim that \refco{main-sch} is equivalent to saying that $\calS(X)=\calJ(G)_U$ where the latter means $U$-coinvariants with respect to the right action of $U$ on $\calJ(G)$ (note that $\calJ(G)$ is a $G$-bimodule, since it contains $\calH(G)$ as a subalgebra). Indeed, let us assume \refco{main-sch}. Then we can define a map $\zeta: \calJ(G)\to \calS(X)$ by sending every $f\in \calJ(G)$ to $\int_U f(gu)du$. The fact that the action of $\calJ(G)$ on $\calS_c(X)$ is well-defined guarantees that this integral is convergent; in fact we have $\zeta(f)=f\star \del_{K/K\cap U}$ for a sufficiently small open compact subgroup $K$ of $G$ (here $\del_{K/K\cap U}$ denotes the multiple of the characteristic function of $K/K\cap U\subset G/U$ normalized by the condition that its integral over $G/U$ is equal to 1; it is easy to see that the result is independent of the choice of $K$ if we require that $f$ is $K$-invariant). Also \refco{main-sch} guarantees that $f\star \del_{K/K\cap U}\in \calS(X)$. It is clear that  $\zeta$ factorizes through $\calJ(G)_U$ and the resulting map $\calJ(G)_U\to \calS(X)$ is injective. On the other hand,  the restriction of $\zeta$ to $\calH(G)$ defines a surjective map $\calH(G)\to \calS_c(X)$. Hence, by definition we $\zeta$ is a surjective map from $\calJ(G)$ to $\calJ(G)\cdot \calS_c(X)=\calS(X)$. Thus we have proved that \refco{main-sch} implies that $\calS(X)=\calJ(G)_U$. The converse statement obvious, since as a $\calJ(G)$-module the space $\calJ(G)_U$ is clearly generated by $\calH(G)_U=\calS_c(X)$.

We can prove the following weaker version of \refco{main-sch}:
\th{main-sch'}
We have $(\calS(X))_{\theta}=(\calJ(G)\cdot \calS_c(X))_{\theta}$
where the character $\theta:\bfT(\calO)\to S^1$ is either trivial or
if the composition of $\theta$ with any coroot is non-trivial.
\eth

\prf
Let $\calS_c=\calS_c(X), \calS=\calS(X),\calS'=\calJ(G)\cdot \calS_c(X)$.
We want to show the equality $\calS_{\theta}=\calS'_{\theta}$ for $\theta$ as above.

\medskip
\noindent
{\em Step 1.}  Let us first show that $\calS^I=(\calS')^I$. As before, let $G^{\vee}$ denote the Langlands dual group of $\bfG$ over $\CC$ and let $\calB$ be its flag variety. According to \cite{Cells2} the algebra $J$ decomposes as a direct sum of subalgebras $J_u$ numbered by unipotent elements $u\in G^{\vee}$ up to conjugacy. We denote by $J_0$ the summand corresponding to the unit conjugacy class. We claim the action on $J$ on $L^2(X)^I$ factorizes through the projection on $J_0$. This is obvious since $L^2(X)^I$ is a torsion-free module over the center of $\calH(G,I)$ and every $J_u$ with non-trivial $u$ is annihilated by a non-zero ideal of the center.

According to \cite{xi-1} the algebra $J_0$ is naturally isomorphic to the $K_{G^{\vee}}(\calB\times \calB)$ (here $K_{G^{\vee}}(?)$ stands for the
complexified Grothendieck group of $G^{\vee}$-equivariant coherent sheaves on ?). On the other hand, let $\calK=K_{T^{\vee}\times \CC^*}(\calB)$ and let $\calK_q$ be its specialization at $v=q^{1/2}$ where $K_{\CC^*}(pt)=\CC[v,v^{-1}]$ (the action of $\CC^*$ on $\calB$ trivial). According to Section 5 of \cite{BK} the space $\calK$ has a natural action of $\calH_{\aff}$ and hence $\calK_q$ has a natural action of $\calH(G,I)$. Moreover, we have an isomorphism $\calS^I\simeq \calK_q$ which identifies
$\calS_c^I$ with the submodule generated by the sky-scraper $\kap$ at some $T^{\vee}$-invariant point $e\in\calB^{\vee}$.
On the other hand $J_0\ten \calA=K_{G^{\vee}\times \CC^*}(\calB\times\calB)$ clearly acts on $\calK$ (and this action is compatible with the $\calH_{\aff}$-action with respect to the homomorphism $\calH_{\aff}\to J_0$ -- this is proved in \cite{xi-1}); moreover the action of $K_{G^{\vee}\times \CC^*}(\calB\times\calB)$ on $\kap$ defines an isomorphism $K_{G^{\vee}\times \CC^*}(\calB\times\calB)\simeq \calK$. Hence the same is true after specialization to $v=q^{1/2}$.
We see that $J_0$ acts on $\calK_q=\calS^I$ and the latter is generated as a module by an element of $\calS_c^I$. This implies the equality $\calS^I=(\calS')^I$.

\medskip
\noindent
{\em Step 2.} Let $\calS_0,\calS'_0$ denote the $G$-module of coinvariants of $\calS$ (resp. $\calS'$) with respect to $\bfT(\calO)$. Then both are subrepresentations of $C^{\infty}(X)_0$. Recall that if a $G$-module $V$ is generated by $I$-fixed vectors then two $G$-submodules $W_1$ and $W_2$
of $V$ coincide if and only if $W_1^I=W_2^I$. Applying this to $C^{\infty}(X)_0, W_1=\calS_0,W_2=\calS'_0$ and using Step 1 we get the equality $\calS_0=\calS'_0$.

\medskip
\noindent
{\em Step 3.}
Let $\theta:\bfT(\calO)\to S^1$ be a character. For a $T$-module $V$ let $V_{\theta}$ denote the corresponding space of $(\bfT(\calO),\theta)$-coinvariants. Let us prove that $\calS_{\theta}=\calS'_{\theta}$ assuming that the composition of $\theta$ with any coroot is a non-trivial character of $\calO^*$. We shall refer to such $\theta$ as "regular".

In this case it is obvious from the definition of $\Phi_w$ that for any simple coroot $\alp$ the operator $\Phi_{s_{\alp}}$ defines an isomorphism between  $\calS_{c,\theta}$ and $\calS_{c,s_{\alp}(\theta)}$. Indeed, any $(\bfT(\calO),\theta)$-equivariant function automatically vanishes on $\oX^{\alp}\backslash X$ (recall the notation of \refss{dig}) and the same is true for $\theta$ replaced with $s_{\alp}(\theta)$.  Since the notion of regularity is $W$-invariant it follows that for any $w\in W$ the operator $\Phi_w$ defines an isomorphism
 between $\calS_{c,\theta}$ and $\calS_{c,w(\theta)}$. On the other hand, we claim that $\calS'_{\theta}=\calS_{c,\theta}$.
For this it is enough to prove that $\calS_{c,\theta}$ is $\calJ(G)$-invariant. This would follow if we knew that for any character $\chi:T\to \CC^*$ such that $\chi|_{\bfT(\calO)}=\theta$ the action of $\calH(G)$ on the space $\calS_{c,\chi}$ of $(T,\chi)$-coinvariants on $\calS_c$ extends to $\calJ(G)$. For any $\chi$ as above we can find an element $w\in W$ such that $w(\chi)$ is non-negative. Hence by definition $\calJ(G)$ acts on $\calS_{c,w(\chi)}=i_{GB}(w(\chi))\simeq i_{GB}(\chi)$.


\epr
\ssec{}{The parabolic case}Let $P$ be a parabolic subgroup of $G$ with a Levi subgroup $M$ and unipotent radical $U_P$.
Let $X_P=G/U_P$. This space has a natural $G\times M$ action. Therefore the space $\calS_c(X)$ of locally constant compactly supported functions on $X_P$ becomes a $G\times M$ module; for convenience we are going to twist the $M$ action by the square root of the absolute value of the determinant of the $M$-action on the Lie algebra $\gru_P$ of $U_P$.

As before, we can define the space $L^2(X_P)$. For the same reason as before it has an action of $\calC(G)$. We now define $\calS(X_P):=\calJ(G)\cdot \calS_c(X_P)\subset L^2(X_P)$. Equivalently $\calS(X_P)=\calJ(G)_{U_P}$.
\conj{parabolic}
Let $P$ and $Q$ be two associate parabolics, i.e. two parabolics with the same Levi subgroup $M$. Then there exists a $G\times M$-equivariant unitary isomorphism $\Phi_{P,Q}:L^2(X_P)\widetilde{\to} L^2(X_Q)$ whose restriction to $\calS(X_P)$ defines an isomorphism between $\calS(X_P)$ and $\calS(X_Q)$.
\econj

\sec{}{Appendix: an $SL(2)$-example}

\ssec{}{}The purpose of this appendix is to show how the algebra $J$ gets realized inside the locally constant functions on $G$ for the case when $G=SL(2,F)$. Let $St$ denote the Steinberg representation of $G$. Let also $\calS$ denote the space of locally constant compactly supported functions on $F^2$. This space has an action of $G$ and a commuting action of the Fourier transform $\Phi$ and of the torus $T=F^{\x}$. Similarly, we let $\calS_c$ denote the space of locally constant compactly supported functions on $F^2\backslash \{ 0\}$.  It follows from the above results that
\eq{J-sl2}
J=\End(St^I)\oplus J_0,\quad \text{where $J_0=\End_{\Phi,F^{\x}}(\calS^I)$}
\end{equation}
The embedding $\calH(G,I)\hookrightarrow J$ is given by the action of $\calH(G,I)$ on $St^I$ and on $\calS^I$.
Since the first summand in \refe{J-sl2} is a one-dimensional subspace that does not lie in $\calH(G,I)$, it follows that the projection
of $\calH(G,I)$ to $J_0$ is an embedding. Moreover, it is clear that the codimension of $\calH(G,I)$ in $J_0$ is at most 1. Indeed, an element of $J_0$ comes from an element of $\calH(G,I)$ if and only if it sends $\calS_c^I$ to $\calS_c^I$. The quotient $\calS/\calS_c$ is naturally isomorphic to $\CC$ (the map is given by evaluating a function at $0$) with $F^{\x}$-action given by the character $x\mapsto q^{|x|}$ (recall the convention about the torus action from \refss{dig}). Denote this character by $\chi$ and let $\calS_{\chi}$ (resp. $(\calS_c)_{\chi}$) denote the space of $(F^{\times},\chi)$-coinvariants on $\calS$ (resp. on $\calS_c$). Let $U$ denote the image of $(\calS_c)_{\chi}$ in $\calS_{\chi}$. Then it is clear that the quotient $J_0/\calH(G,I)$ embeds into $\Hom_T(U^I,\CC)$. However, $U$ is actually isomorphic to $St$, hence $\dim U^I=1$.
This shows that $\dim J_0/\calH(G,I)$ is either $0$ or $1$ and in the latter case it is isomorphic to $\CC\otimes St^I$ as an $\calH(G,I)$-bimodule.

We would like to show that $J_0/\calH(G,I)$ is indeed one-dimensional. Note that this will imply that $\dim J/\calH(G,I)=2$.
For this it is enough to find one $I\times I$-invariant function $f$ on $G$ such that

(a) $f$ does not have compact support but belongs to $\calC(G,I)$,

(b) The action of $f$ on $\calS^I$ is well-defined and it is given by a non-zero operator.

\noindent
We are going to find such a function explicitly. In fact, the function $f$ will be $K\times I$-invariant, where $K=SL(2,\calO)$.
We begin by some generalities.
\ssec{}{General remarks and volume calculation}
Let $G$ be a totally disconnected group, $x\in G$, $K_1,K_2$ - open compact subgroups of $G$; set $X=K_1 x K_2$.
Let $V$ be a smooth representation of $G$, and $v\in V$. Fix a Haar measure on $G$. Then we have an element $e_X$ in the Hecke algebra of $G$.
By definition, we have
$$
e_X(v)=\int\limits_{g\in X} g(v) dg.
$$
This formula can be rewritten in the following way, which will be crucial in the future:
\eq{formula'}
e_X(v)=\frac{\vol(X)}{\vol(K_1)\cdot \vol(K_2)}\int\limits_{g_1\in K_1,g_2\in K_2} (g_1 x g_2)(v) dg_1 dg_2.
\end{equation}
In particular, let us assume that the Haar measure is chosen is such a way that $\vol(K_2)=1$ and that $v$ is $K_2$-invariant.
Then \refe{formula'} takes the form
\eq{formula}
e_X(v)=\frac{\vol(X)}{\vol(K_1)}\int\limits_{g\in K_1} g (x(v)) dg.
\end{equation}

We now want to apply this formula to the case when $G=SL(2,F), K=K_1=SL(2,\calO), K_2=I$ (so, in particular, we fix the Haar measure so that the volume of $I$ is 1).
Let $t$ be a uniformizer of $F$ and let $x_n=\begin{pmatrix} t^n & 0\\ 0 & t^{-n}\end{pmatrix}$ where $n\in \ZZ$ and let $X_n=K x_n I$.
\lem{}
We have
$$
\frac{\vol(X_n)}{\vol(K)}=
\begin{cases}
q^{2n-1}\ \text{if $n\geq 0$}\\
q^{-2n}\ \text{if $n < 0$}
\end{cases}
$$
\elem
\prf
Let $H_n\subset K$ be the subgroup of all $h\in K$ such that $x_n^{-1} h x_n\in I$. Then $\vol(X_n)=\# (K/H_n)$. Recall that $\vol(K)=\# (K/I)=(q+1)$.

If $h=\begin{pmatrix} a & b\\ c & d\end{pmatrix}$ then $x_n^{-1} h x_n=\begin{pmatrix} a & t^{-2n}b\\ t^{2n}c & d\end{pmatrix}$. Hence,
if $n > 0$, then $x_n^{-1} h x_n\in I$ iff $b\in t^{2n}\calO$ and the cardinality of $K/H_n$ is $q^{2n-1}(q+1)$. If $n\leq 0$ then the condition is $c\in t^{-2n+1}\calO$ and the cardinality of $K/H_n$ is $q^{-2n}(q+1)$.
\epr
We now pass to the function $f$.
\prop{app-main}
Let
\eq{gamn}
\gam_n=
\begin{cases}
q^{2n}\ \text{if $n\leq 0$}\\
-q^{-2n+1}\ \text{if $n > 0$}
\end{cases}
\end{equation}
Let $f=\sum_{n\in \ZZ} \gam_n \cdot \chi_{X_n}$ where $\chi_{X_n}$ denotes the characteristic function of $X_n$. Then $f$ belongs to $J_0$ (and it obviously does not belong to $\calH(G,I)$).
\eprop
\cor{}
As a bimodule over $\calH(G,I)$ the quotient $J/\calH(G,I)$ is isomorphic to $V\otimes St^I$ where $V$ is a 2-dimensional representation of $\calH(G,I)$ which is a non-trivial extension of $\CC$ (trivial representation) by $St^I$.
\ecor
\prf
We have already explained above that \refp{app-main} implies that $\dim J/\calH(G,I)=2$. Moreover, it follows from \refe{J-sl2} that this quotient contains $St^I\otimes St^I$ as an $\calH(G,I)$-submodule and we have explained at the end of \refss{} that the quotient of $J/\calH(G,I)$ by $St^I\otimes St^I$ is isomorphic to $\CC\otimes St^I$ as an $\calH(G,I)$-bimodule. Hence we see that $J/\calH(G,I)$ is isomorphic to $V\otimes St^I$ where $V$ is a 2-dimensional representation of $\calH(G,I)$ which is an extension of $\CC$ by $St^I$. If that extension were trivial, then $f$ would have be $G$-invariant on the left modulo functions with compact support. This would imply that $\gam_n$ is constant for $|n|>>0$, which is obviously not the case (and it is also clear that no element of $\calC(G)$ could satisfy this).
\epr
Let us now turn to the proof of \refp{app-main}.
It is easy to see that $f$ belongs to $\calC(G)$. Also since $f$ is $K$-invariant on the left, its projection to the first summand of \refe{J-sl2} is 0. Hence to show that it belongs to $J_0$ it is enough to show that its action on $\calS^I$ is well-defined.
For this, we are going to take $f$ of the form $\sum \gam_n \chi_{X_n}$ with arbitrary coefficients $\gam_n$ and see what conditions we need to impose of the coefficients so the action of $f$ on $\calS^I$ is well-defined. To check the latter, it is enough to verify that
$f\star \chi_{\calO\oplus\calO}$ and $f\star \chi_{\calO\oplus t\calO}$ are well-defined.
\subsubsection{The case of $\chi_{\calO\oplus\calO}$}
Let us first compute the action of $f$ on $\chi_{\calO\oplus \calO}$. It is enough to compute the value of the action of $\chi_{X_n}$ on $\chi_{\calO\oplus \calO}$ at $(t^{-r},0)$ where $r\in \ZZ$.
It is equal to $\frac{\vol(X_n)}{\vol(K)}\cdot \vol(K_{n,r})$ where $K_{n,r}=\{ h\in K|\ h(t^{-r},0)\in t^n\calO\oplus t^{-n}\calO\}$.

Assume first that $n>0$.
Then we have the following cases:
\begin{enumerate}
\item
$r> n$. Then the result is $0$.
\item
$r\leq -n$. Then $K_{n,r}=K$ and the result is $(q+1)q^{2n-1}$.

\item
$-n < r\leq n$. In this case if $h=\begin{pmatrix} a & b\\ c & d\end{pmatrix}$ then $h\in K_{n,r}$ iff $a\in t^{n+r}\calO$ and thus
$\vol(K_{n,r})=q^{-n-r+1}$. Hence the value is $q^{n-r}$.
\end{enumerate}

On the other hand, assume that $n=-m, m>0$. Then we have the following cases:
\begin{enumerate}
\item
$r>m$ -- the value is 0

\item
$r\leq -m$. Then again $K_{n,r}=K$ and we get $(q+1)q^{2m}$.

\item
$-m< r\leq m$. Then $h\in K_{n,r}$ iff $c\in t^{r+m}\calO$ and its volume is $q^{-r-m+1}$. Hence the value is $q^{m-r+1}$.
\end{enumerate}
If $n=0$ then we just get $(q+1)\chi_{\calO\oplus\calO}$.
Altogether, the value of $f\star \chi_{\calO\oplus\calO}$ at $(t^{-r},0)$ for $r>0$ is
$$
\sum\limits_{n\geq r} \gam_n q^{n-r} +\sum\limits_{m\geq r} \gam_{-m} q^{-m-r+1}.
$$
Hence if we want this to be $0$ (which we definitely want at least for $r$ large enough)
we need
\eq{rel1}
\gam_r + q \gam_{-r} =0.
\end{equation}

\subsubsection{The case of $\chi_{\calO\oplus t\calO}$}
In this case we define $K_{n,r}'=\{ h\in K|\ h(t^{-r},0)\in t^n\calO\oplus t^{-n+1}\calO\}$ and we need to repeat the above calculation with $K_{n,r}$ replaced by $K_{n,r}'$.

Let us first assume that $n>0$. Then we have the following cases:
\begin{enumerate}
\item
$r< n-1$ -- get $0$
\item
$r\leq -n$ -- get $(q+1)\frac{\vol(X_n)}{\vol(K)}\gam_n=(q+1)q^{2n-1}\gam_n$

\item  $\-n < r \leq n-1$

In this case $K_{n,r}'$ is given by the condition $a\in t^{n+r}\calO$ and the value is $q^{-n-r+1}q^{2n-1}\gam_n=q^{n-r}\gam_n$.

\end{enumerate}

On the other hand, assume that $n\leq 0$ and set $n=-m$. Then we have the following cases:
\begin{enumerate}
\item $r>m$ -- get $0$

\item $r\leq -m-1$ -- get $(q+1)q^{2m}\gam_{-m}$.

\item
$-m\leq r\leq m$. Then $K_{n,r}'$ is given by the condition $c\in t^{r+m+1}\calO$ and the volume of $K_{n,r}'$ is $q^{-r-m}$.
Thus we get $q^{-r-m}q^{2m}\gam_{-m}=q^{m-r}\gam_{-m}$.
\end{enumerate}

\medskip
\noindent
Thus for $r\geq 0$ we get
that the value of $f\star \chi_{\calO\oplus t\calO}$ at $(t^{-r},0)$ is equal to
$$
\sum\limits_{n\geq r+1} q^{n-r}\gam_n +\sum\limits_{m\geq r} q^{m-r}\gam_{-m}.
$$
This going to be $0$ for all $r\geq 0$ if
\eq{rel2}
q\gam_{r+1}+ \gam_{-r}=0.
\end{equation}

\ssec{}{The function $f$}
Altogether it is now clear that if $\gam_n$ is given by \refe{gamn} then $f$ satisfies both \refe{rel1} and \refe{rel2}.

Let us now compute the action of the resulting function $f$ on both $\chi_{\calO\oplus \calO}$ and $\chi_{\calO\oplus t\calO}$.
In the first case, we already know that $f\star \chi_{\calO\oplus\calO}(t^{-r},0)=0$ if $r>0$. On the other hand, for $r\leq 0$ we have
$$
\begin{aligned}
&f\star \chi_{\calO\oplus \calO}(t^{-r},0)=\\
&(q+1)+\sum\limits_{0<n\leq -r} (q+1) q^{2n-1}(-q^{-2n+1})+\sum\limits_{n>-r} q^{n-r}(-q^{-2n+1})+\\ &\sum\limits_{0<m\leq -r} (q+1)q^{2m}q^{-2m}+
\sum\limits_{m>-r} q^{m-r+1} q^{-2m}=q+1,
\end{aligned}
$$
which shows that $f\star \chi_{\calO\oplus\calO}=(q+1)\chi_{\calO\oplus\calO}$.
Similar calculation shows that $f\star \chi_{\calO\oplus t\calO}=0$.


\begin{thebibliography}{ddddd}
\bibitem{Ber}
J.~Bernstein, {\em  Draft of : Representations of p-adic groups}, Fall 1992. Lectures by Joseph Bernstein,
Written by Karl E. Rumelhart.

\bibitem{BK}
A.~Braverman and D.~Kazhdan,
{\em On the Schwartz space of the basic affine space},
Selecta Math. (N.S.) {\bf 5} (1999), no. 1, 1-28.

\bibitem{KL}
D.~Kazhdan and G.~Lusztig,
{\em Proof of the Deligne-Langlands conjecture
for Hecke algebras}, Invent. math. {\bf 87}, 153-215 (1987)


\bibitem{Cells2}
G.~Lusztig,  {\em Cells in affine Weyl groups. II.} J. Algebra {\bf 109} (1987), no. 2, 536-548.

\bibitem{Cells3}
G.~Lusztig,  {\em Cells in affine Weyl groups. III.}
 J. Fac. Sci. Univ. Tokyo Sect. IA Math. {\bf 34} (1987), no. 2, 223-243.

\bibitem{Cells4}
G.~Lusztig,  {\em Cells in affine Weyl groups. IV.}
 J. Fac. Sci. Univ. Tokyo Sect. IA Math. {\bf 36} (1989), no. 2, 297-328.

\bibitem{sil}
A.~J.~Silberger,
{\em Introduction to Harmonic Analysis on Reductive P-Adic Groups}, Based on Lectures by Harish-Chandra at The Institute for Advanced Study, 1971-73. Princeton University Press, 1979.

\bibitem{Wal}
J.-L. Waldspurger,
{\em La formule de Plancherel pour les groupes p-adiques
(d'apr\`es Harish-Chandra)}. J. Inst. Math. Jussieu, 2(2):235-333,
2003

\bibitem{xi}
N.~Xi, {\em Representations of affine Hecke algebras and based rings of affine Weyl groups.} J. Amer. Math. Soc. {\bf 20} (2007), no. 1, 211-217.

\bibitem{xi-1}
N.~Xi, {\em The based ring of the lowest two-sided cell of an affine Weyl group, III}, arXiv:1506.00476.

\end{thebibliography}
\end{document}